\documentclass[reqno]{amsart}

\theoremstyle{plain}
  \newtheorem{theo}{Theorem}
  \newtheorem{coro}{Corollary}
  \newtheorem{prop}{Proposition}
  \newtheorem{lemm}{Lemma}
  \newtheorem*{quest*}{Question}

\theoremstyle{remark}
  \newtheorem*{rema*}{Remark}

\newcommand\itemref[1]{(\ref{#1})}

\newcommand\resp{resp.\ }
\newcommand\ie{i.e.\ }
\newcommand\cf{cf.\ }

\newcommand\Q{\mathbb Q}
\newcommand\C{\mathbb C}
\newcommand\R{\mathbb R}
\newcommand\Z{\mathbb Z}

\newcommand\boe{\mathfrak b}
\newcommand\goe{\mathfrak g}
\newcommand\hoe{\mathfrak h}
\newcommand\sloe{\mathfrak{sl}}

\newcommand\Vhoe{{\mathcal V_\hoe}}
\newcommand\Vboe{{\mathcal V_\boe}}
\newcommand\Vgoe{{\mathcal V_\goe}}

\begin{document}

\title{A remark on the c--splitting conjecture}

\author{Stefan Haller}

\address{Department of Mathematics, University of Vienna, 
         Strudlhofgasse 4, A-1090 Vienna, Austria.}

\email{Stefan.Haller@univie.ac.at}

\thanks{The author is supported by the `Fonds zur F\"orderung der 
        wissenschaftlichen Forschung' (Austrian Science Fund),
        project number {\tt P14195-MAT}}



\begin{abstract} 
Let $M$ be a closed symplectic manifold and suppose
$M\to P\to B$ is a Hamiltonian fibration. Lalonde and 
McDuff raised the question whether one always has
$H^*(P;\Q)=H^*(M;\Q)\otimes H^*(B;\Q)$ as vector spaces. This is known as the
c--splitting conjecture. They showed, that
this indeed holds whenever the base is a sphere. Using their theorem we will 
prove the c--splitting conjecture for arbitrary base $B$ and
fibers $M$ which satisfy a weakening of the Hard Lefschetz condition.
\end{abstract}

\maketitle

\section{Introduction and statement of the result}

Let $M$ be a closed symplectic manifold and consider the group of 
diffeomorphisms preserving the symplectic structure. As a normal subgroup we
find the group of Hamiltonian diffeomorphisms. It consists of all
diffeomorphisms which are integrals of time dependent Hamiltonian vector
fields. Particularly it is connected. 

The group of Hamiltonian diffeomorphisms has many intriguing properties.
For instance every Hamiltonian diffeomorphism is known to have lots of fixed
points, \ie more than the Lefschetz fixed point theorem guaranties
for mappings homotopic to the identity. More precisely, if all
the fixed points are non-degenerate there have to be at least $\sum b_i$
many, where $b_i$ denotes the $i$-th Betti number of $M$. In contrast the
Lefschetz theorem just gives an estimate by $\sum(-1)^ib_i$.
This is a deep theorem with major contributions from Floer, Hofer, Zehnder,
Salamon, Fukaya, Ono, Liu and Tian -- just to name a few.

This note is about another property the group of Hamiltonian
diffeomorphisms seems to have. Recall that a 
Hamiltonian fibration is a fiber bundle $M\to P\to B$ with
typical fiber a closed symplectic manifold whose structure group is reduced to the
group of Hamiltonian diffeomorphisms. For a Hamiltonian fibration one can
show, that the cohomology class of the symplectic structure lies in the image
of $H^*(P)\to H^*(M)$. Conversely, if the structure group of a fiber bundle
$M\to P\to B$ can be reduced 
to the connected component of the group of symplectic diffeomorphisms this
condition in turn implies that the structure group can actually be reduced 
to the group of Hamiltonian
diffeomorphisms. All this can be found in \cite{LM02}.

One says a Hamiltonian fibration
c--splits (short for cohomologically splits) if the cohomology of the total
space satisfies 
\begin{equation}\label{csplit}
H^*(P;\Q)=H^*(M;\Q)\otimes H^*(B;\Q)
\end{equation}
as $H^*(B)$--modules.\footnote{Usually one just asks \eqref{csplit} to hold as
vector spaces. In view of the Leray--Hirsch theorem this is equivalent to
our definition as long as $H^*(B)$ is finite dimensional. However, for
infinite dimensional $H^*(B)$ our condition seems to be more adequate.}
In other words, cohomologically -- disregarding the ring structure --
the fibration looks like a product. It makes no difference if we take
cohomology with reel coefficients. From now on all cohomology groups
are understood to be with coefficients in $\R$, and we will omit them in our
notation.

In \cite{LM02}, Lalonde and McDuff raised the following

\begin{quest*}[Lalonde and McDuff]
Does every Hamiltonian fibration c--split?
\end{quest*}

The affirmative answer to their question is known as the c--splitting 
conjecture. It has been proved to be true
in many circumstances, yet the general case is still a mystery.
In their paper \cite{LM02} Lalonde and McDuff proved, that the c--splitting
conjecture indeed holds whenever the base is a sphere or a $3$--dimensional
CW--complex. The difficult part is the case $B=S^2$, which requires hard
analytic tools, see \cite{LMP99} and \cite{M00}.

Using Lalonde and McDuff's theorem K\c edra derived, that the c--splitting 
conjecture holds for $4$--dimensional fibers, simply connected
$6$--dimensional fibers and simply connected spherically generated 
fibers, see \cite{Ka} and \cite{Kb}.
Employing parameterized Gromov--Witten invariants he also showed,
that the c--splitting conjecture is true whenever the fiber is $\C P^5$ blown
up along Thurston's $4$--dimensional nil-manifold, see \cite{Ka}.

Another situation when the c--splitting conjecture is known to hold, is when 
the structure group reduces to a compact subgroup of the Hamiltonian
diffeomorphisms, see \cite{LM02} and \cite{AB84}.

Finally the c--splitting conjecture holds for fibers which
satisfy the Hard Lefschetz Theorem.
This was already observed by Blanchard, see \cite{B56}.

The purpose of this note is to establish the c--splitting conjecture for a
class of fibers, which satisfy a weakening of the Lefschetz condition.

\begin{theo}\label{main-theo}
Suppose $(M,\omega)$ is a closed symplectic manifold of dimension $2n$, such 
that the Lefschetz type mappings
$$
[\omega]^k:H^{n+1-k}(M)\to H^{n+1+k}(M)
$$
are onto for all $k\geq 0$. 
Then every Hamiltonian fibration $M\to P\to B$ c--splits. 
\end{theo}

Let us remark, that the main ingredient in our proof is Lalonde and McDuff's
theorem which tells, that every Hamiltonian fibration over a $3$--dimensional
CW--complex c--splits. We then apply methods which are in essence the same
Blanchard used to proof the c--splitting for Lefschetz fibers.
However, we hope our approach is easier to use and more conceptual.

It is easy to see, that a fiber bundle $M\to P\to B$
c--splits iff $H^*(P)\to H^*(M)$ is onto or equivalently iff 
$H_*(M)\to H_*(P)$ is injective. Also the bundle will c--split iff the
Leray--Serre spectral sequence collapses at the $E^2$--term, \ie its
differentials $\partial_k:E^k\to E^k$ vanish for all $k\geq 2$. 

Essentially it suffices to consider bases $B$ which are finite CW--complexes. 
Indeed, fix a closed symplectic fiber $M$ and suppose every Hamiltonian 
fibration with fiber $M$ and a finite CW--complex as a base c--splits.
From the homological interpretation above it is clear, that
this implies the c--splitting conjecture for arbitrary bases $B$ and fiber
$M$.

Particularly we can look at the universal Hamiltonian fibration. Fix a
closed symplectic manifold $M$ and let $G$ denote the group of Hamiltonian
diffeomorphisms. Let $G\to EG\to BG$ denote the universal $G$ bundle and
consider the associated universal Hamiltonian fibration 
\begin{equation}\label{universal}
M\to EG\times_GM\to BG.
\end{equation}
Whenever $M\to P\to B$ is another Hamiltonian fibration with the
same fiber there is a map 
$f:B\to BG$, such that $P=f^*(EG\times_GM)$. One easily derives, that if the
c--splitting conjecture holds for \eqref{universal} it will hold for all
Hamiltonian fibrations with fiber $M$.

The cohomology of the total space $EG\times_GM$ is 
known as the equivariant cohomology of $M$ with respect to the 
action of $G$. So the bundle \eqref{universal} will c--split if and only if
the equivariant cohomology is a free module over $H^*(B)$.
So the conjecture of Lalonde and McDuff can be reformulated as follows:
For every closed symplectic manifold the equivariant cohomology of $M$
with respect to the action of the Hamiltonian group is a free
$H^*(B)$--module.

Finally let us remark, that the c--splitting property is a geometric rather
than a topological phenomenon. In \cite{LM02} Lalonde and McDuff constructed
a smooth fiber bundle $M\to P\to S^2$ with $6$--dimensional closed
fiber. Its total space admits a class $\alpha\in H^2(P)$ which satisfies 
$0\neq\alpha^3\in H^6(M)$. In some sense this is the cohomological
analogue of a Hamiltonian fibration. However this bundle does not
c--split. 

\section{Canonic filtration of $\boe$-modules}

Let $\goe:=\sloe(2;\R)$ with base $\{e,f,h\}$ and relations $[h,e]=2e$,
$[h,f]=-2f$, $[e,f]=h$. Let $\hoe$ denote the subalgebra spanned by $h$, and
$\boe$ the subalgebra spanned by $\{e,h\}$. Let $\Vhoe$ denote the category of
$\hoe$--modules $V$, which admit a decomposition $V=\bigoplus_{k\in\Z}V^k$
into eigenspaces of $h$, $V^k$ being the eigenspace to the weight $k$,
and only finitely many $V^k$ non-trivial. Moreover let $\Vboe$ \resp $\Vgoe$
denote the category of $\boe$ \resp $\goe$--modules for which the underlying
$\hoe$--module is in $\Vhoe$. Then $e:V^k\to V^{k+2}$ and $f:V^k\to
V^{k-2}$.

In this section we will collect a few basic properties of $\boe$--modules
which we are going to use in the proof of Theorem~\ref{main-theo}. Most importantly
the existence of a canonic filtration for every $V\in\Vboe$. This filtration
was used by Mathieu \cite{M95} when he proved that a symplectic
manifold satisfies the Hard Lefschetz Theorem iff every cohomology class has
a harmonic representative in the sense of Brylinksi, see \cite{B88}.
The proofs for all the statements below are elementary and can be found 
in \cite{H}, see also \cite{M95}.

\begin{lemm}\label{boe-vs-goe-lem}
Suppose $V,W\in\Vgoe$ and $\varphi:V\to W$ a $\boe$--module homomorphism.
Then $\varphi$ is a $\goe$--module homomorphism.
\end{lemm}

For $V\in\Vboe$ we write $V\in\Vgoe$
if the $\boe$--module structure extends to a $\goe$--module structure. The
previous lemma tells, that such a $\goe$--module structure is unique if
it exists.

For $V\in\Vboe$ and $k\in\Z$ we will denote by $V[k]$ the $\boe$--module 
which has $V$ as underlying vector space, the action of $e\in\boe$ is the
same as on $V$ but the $\hoe$--action is shifted by $k$, \ie
$h\cdot v=hv+kv$. Here $h\cdot v$ is supposed to denote the new
$\hoe$--action on $V[k]$, whereas $hv$ denotes the old $\hoe$--action on $V$.

\begin{prop}\label{filtration}
Suppose $V\in\Vboe$. Then there exists a unique filtration 
$\cdots\subseteq V_{m-1}\subseteq V_m\subseteq\cdots$ 
of $V$ with the following properties:
\begin{enumerate}
\item
$V_m=0$ for $m$ sufficiently small.
\item
$V_m=V$ for $m$ sufficiently large.
\item
$V_m\subseteq V$ is a $\boe$--submodule, for all $m\in\Z$.
\item
$(V_m/V_{m-1})[-m]\in\Vgoe$, for all $m\in\Z$. 
\end{enumerate}
\end{prop}

One readily verifies the following 

\begin{lemm}\label{dsp}
Suppose $V,W\in\Vboe$. Then:
\begin{enumerate}
\item\label{filt-dual}
$(V^*)_m=\{\alpha\in V^*:\alpha|_{V_{-m-1}}=0\}$.
\item\label{filt-sum}
$(V\oplus W)_m=V_m\oplus W_m$.
\item\label{filt-prod}
$(V\otimes W)_m=\sum_{m_1+m_2=m}V_{m_1}\otimes W_{m_2}$.
\item\label{filt-shift}
$(V[k])_{m+k}=V_m$.
\end{enumerate}
\end{lemm}

\begin{prop}\label{filtration_preserving}
Suppose $V,W\in\Vboe$ with corresponding filtrations
$V_m$ and $W_m$. Then every $\boe$--module homomorphism $\varphi:V\to W$ 
is filtration preserving, that is
$\varphi(V_m)\subseteq W_m$, for all $m\in\Z$.
\end{prop}

\begin{coro}\label{filt-pres-coro}
Let $V,W\in\Vboe$ with corresponding filtrations $V_m$ and $W_m$. Suppose
$\varphi:V\to W$ is a linear map satisfying $\varphi(ev)=e\varphi(v)$ and
$\varphi(hv+kv)=h\varphi(v)$, for all $v\in V$ and some fixed $k\in\Z$.
Then $\varphi(V_m)\subseteq W_{m+k}$.
\end{coro}

\begin{proof}
The assumption on the map $\varphi:V\to W$ is equivalent to 
$\varphi:V[k]\to W$ being a
$\boe$--module homomorphism. Using Lemma~\ref{dsp}\itemref{filt-shift} and
Proposition~\ref{filtration_preserving} we conclude 
$\varphi(V_m)=\varphi((V[k])_{m+k})\subseteq W_{m+k}$. 
\end{proof}

\begin{prop}\label{bry-prop}
Let $V\in\Vboe$, $m\in\Z$ and let $V=\bigoplus_{k\in\Z}V^k$ denote the
decomposition of $V$ into eigenspaces of $h$. Then the following are
equivalent:
\begin{enumerate}
\item\label{bp-i}
$V=V_m$.
\item\label{bp-ii}
$e^l:V^{m-l}\to V^{m+l}$ is onto for all $l\geq 0$.
\end{enumerate}
\end{prop}

Finally let us remark, that for a finite dimensional $V\in\Vboe$ one can
give explicit formulas for the dimensions of $V^k_m$ in terms of the ranks
of all the mappings $e^i:V^j\to V^{j+2i}$. This can be found in \cite{H} but
we won't make use of it in the sequel.

\section{Proof of Theorem~\ref{main-theo} and examples}\label{proof}

Let $M$ be a topological space and $\alpha\in H^2(M)$. Consider the 
cohomology $H^*(M)$ as a $\boe$--module via
$e\cdot\beta:=\alpha\cup\beta$ for $\beta\in H^*(M)$ and
$h\cdot\beta:=k\beta$ for $\beta\in H^k(M)$. Let $H^*(M)_m$ denote the
corresponding filtration from Proposition~\ref{filtration}.
The next proposition can be expressed most conveniently using the 
associated graded space $\tilde H^*(M)_m:=H^*(M)_m/H^*(M)_{m-1}$.

\begin{prop}[Poincar\'e duality]\label{poincare}
Suppose $M$ is an closed oriented manifold of dimension $n$, 
$\alpha\in H^2(M)$ and $m\in\Z$.
Then the Poincar\'e pairing factors to a non-degenerate bilinear pairing
\begin{equation}\label{poin}
\tilde H^*(M)_m\otimes\tilde H^*(M)_{n-m}\to\R,
\quad
\beta\otimes\gamma\mapsto(\beta\cup\gamma)\cap[M].
\end{equation}
\end{prop}

\begin{proof}
Consider the Poincar\'e duality 
$$
\Phi:H^*(M)\to\bigl(H^*(M)[-n]\bigr)^*,
\quad
\Phi(\beta)(\gamma)=(\beta\cup\gamma)\cap[M]
$$
and the mapping
$$
\Psi:H^*(M)\to H^*(M),
\quad
\Psi(\beta)=(-1)^{k(k+1)/2}\beta
\quad\text{for $\beta\in H^k(M)$.}
$$
One easily checks, that $\Phi\circ\Psi$ is 
a $\boe$-module homomorphism. From Proposition~\ref{filtration_preserving} 
we thus get
$$
\Phi(H^*(M)_m)
=\Phi(\Psi(H^*(M)_m))
\subseteq((H^*(M)[-n])^*)_m.
$$ 
Using Lemma~\ref{dsp} we conclude that \eqref{poin} is
well defined for every $m\in\Z$. 
It follows from ordinary Poincar\'e duality, that it has to be non-degenerate.
\end{proof}

\begin{coro}\label{lef-vs-fil}
Let $M$ be an oriented closed manifold of dimension $n$, $\alpha\in H^2(M)$ 
and $m\in\Z$. Then the following are equivalent:
\begin{enumerate}
\item
$\alpha^k:H^{m-k}(M)\to H^{m+k}(M)$ is onto, for all $k\geq 0$.
\item
$H^*(M)_m=H^*(M)$. 
\item
$H^*(M)_{n-m-1}=0$.
\end{enumerate}
\end{coro}

\begin{proof}
The equivalence of the first two assertions is an application of
Proposition~\ref{bry-prop}. The last two statements are 
equivalent, for we have Proposition~\ref{poincare}.
\end{proof}

We are now in a position to apply the algebraic machinery and prove
Theorem~\ref{main-theo}.
As warm up exercise we give a proof of the c--splitting conjecture for
fibers which satisfy the Hard Lefschetz Theorem. As was already mentioned in
the introduction, this is
an old theorem due to Blanchard, see \cite{B56}. Recall, that a symplectic
manifold $M$ of dimension $2n$ is said to satisfy the Hard Lefschetz Theorem 
if the Lefschetz maps
$$
[\omega]^k:H^{n-k}(M)\to H^{n+k}(M)
$$
are onto, for all $k\geq 0$. Corollary~\ref{lef-vs-fil} tells us, that for a 
closed oriented $M$ this condition is
equivalent to $H^*(M)_n=H^*(M)$ and $H^*(M)_{n-1}=0$, where we
consider $H^*(M)$ with the $\boe$--module structure induced from
$[\omega]\in H^2(M)$ as described above.

\begin{theo}[Blanchard]\label{blan}
Suppose $(M,\omega)$ is a closed symplectic manifold of dimension $2n$
which satisfies the Hard Lefschetz Theorem. Then every Hamiltonian fibration 
$M\to P\to B$ c--splits.
\end{theo}

\begin{proof}
Consider the Leray--Serre spectral sequence of $P$. We 
consider its $E^2$--term $E^2=H^*(M)\otimes H^*(B)$ as $\boe$--module as
follows. Equip $H^*(M)$ with the $\boe$--module structure induced from 
$[\omega]\in H^2(M)$, $H^*(B)$ with the trivial
$\boe$--module structure and put the tensor product structure on $E^2$.

Since $(M,\omega)$ is supposed to satisfy the Hard Lefschetz Theorem
Corollary~\ref{lef-vs-fil}
yields $H^*(M)_n=H^*(M)$ and $H^*(M)_{n-1}=0$.
Applying Lemma~\ref{dsp}\itemref{filt-prod} we thus have:
\begin{equation}\label{e2-filt}
E^2_m=
\begin{cases}
E^2 & m\geq n
\\
0 & m<n
\end{cases}
\end{equation}

Since the fibration is Hamiltonian $[\omega]\in H^2(M)$ is the restriction
of a class in $H^2(P)$.
So $\partial_2(\omega)=0$, where $\partial_2:E^2\to E^2$ denotes the
differential of the $E^2$--term. Since $\partial_2$ is a derivation we
obtain
$\partial_2(\omega\cup\alpha)=\omega\cup\partial_2(\alpha)$ for all
$\alpha\in E^2$. In other words $\partial_2(e\alpha)=e\partial_2(\alpha)$.
Moreover
$$
\partial_2:H^k(M)\otimes H^*(B)\to H^{k-1}(M)\otimes H^*(B),
$$ 
which is the same as saying
$\partial_2((h-1)\alpha)=h\partial_2(\alpha)$.
From Corollary~\ref{filt-pres-coro} we
thus conclude $\partial_2(E^2_m)\subseteq E^2_{m-1}$.
In view of \eqref{e2-filt} this implies $\partial_2=0$.

So we have $E^3=E^2=H^*(M)\otimes H^*(B)$. We equip
$E^3$ with the $\boe$--module structure we used on $E^2$. 
The same arguments as above imply, that the differential
$\partial_3:E^3\to E^3$ satisfies
$\partial_3(E^3_m)\subseteq E^3_{m-2}$
and thus $\partial_3=0$.
Similarly one goes on and shows $\partial_k=0$ for all $k\geq 2$.
\end{proof}

The proof of Theorem~\ref{main-theo} is similar, but will make use of the 
following deep theorem due to Lalonde and McDuff, see \cite{LM02}.

\begin{theo}[Lalonde and McDuff]\label{csplits2}
Every Hamiltonian fibration with base $S^n$ c--splits.
Moreover every Hamiltonian fibration over a $3$--dimensional 
CW--complex c--splits.
\end{theo}

The difficult part here, is to show that this is true for $S^2$. They
manage to do this using Gromov--Witten invariants and Seidel's representation 
of the fundamental group of the Hamiltonian diffeomorphisms
on the quantum cohomology ring of $M$, see \cite{LMP99} and \cite{M00}. 
The other cases are deduced using topological methods.

\begin{proof}[Proof of Theorem~\ref{main-theo}]
Again we consider the Leray--Serre spectral sequence of the fibration.
Theorem~\ref{csplits2} immediately implies, that 
$$
E^2=E^3=E^4=H^*(M)\otimes H^*(B).
$$
We endow $E^4$ with the $\boe$--module structure we used in the proof of
Theorem~\ref{blan}. Via Corollary~\ref{lef-vs-fil} we see,
that the condition on $M$ is equivalent to
$H^*(M)_{n-2}=0$ and $H^*(M)_{n+1}=H^*(M)$. As before we conclude:
\begin{equation}\label{e4filt}
E^4_m=
\begin{cases}
E^4 & m\geq n+1
\\
0 & m<n-1
\end{cases}
\end{equation}
Again, since the fibration is Hamiltonian we
have $\partial_4(\omega)=0$, hence $\partial_4$ commutes with the action of
$e\in\boe$. Using the fact, that
$$
\partial_4:H^k(M)\otimes H^*(B)\to H^{k-3}(M)\otimes H^*(B)
$$
we get $\partial_4(E^4_m)\subseteq E^4_{m-3}$ from
Corollary~\ref{filt-pres-coro}. In view of \eqref{e4filt} we thus must have 
$\partial_4=0$. Similarly one shows $\partial_k=0$, for all $k\geq 4$.
\end{proof}

\begin{rema*}
Let $\varphi:V\to W$ be a $\boe$--module homomorphism. If we know $V_m\neq
0$ implies $W_m=0$ for all $m\in\Z$ we can conclude that $\varphi$ vanishes. 
In essence, that is what we used in the proofs above. 
However, even if $V_m\neq 0$ and $W_m\neq 0$ our method still gives some
information about $\varphi$. 
Since $\varphi$ is filtration preserving it induces $\goe$--module
homomorphisms
$$
\varphi_m:(V_m/V_{m-1})[-m]\to (W_m/W_{m-1})[-m].
$$ 
Now Schur's lemma gives strong restrictions on such mappings. For
instance, if every highest weight which occurs in the left hand side
$\goe$-representation does not 
occur on right hand side we can conclude $\varphi_m=0$.
\end{rema*}

Let us close this note with some examples of symplectic manifolds
satisfying the condition of Theorem~\ref{main-theo}.

For a $4$--dimensional closed symplectic manifold this condition
is trivially satisfied. So every Hamiltonian fibration with
$4$--dimensional fiber c--splits. This was already observed by  K\c edra
as a consequence of Theorem~\ref{csplits2}, see \cite{LM02}.

A $6$--dimensional closed symplectic manifold satisfies the assumption
of Theorem~\ref{main-theo} iff the mapping $[\omega]:H^3(M)\to H^5(M)$
is onto. Via Poincar\'e duality this is equivalent to 
$[\omega]:H^1(M)\to H^3(M)$ being injective. 
Particularly this applies for simply connected $M$, a fact also observed 
by K\c edra.

Salamon gave a
classification of all $6$--dimensional nilpotent Lie algebras, see \cite{S01}.
Since there is exactly one nil-manifold to every nilpotent Lie algebra this
is a classification of all $6$--dimensional nil-manifolds. Many of
them admit symplectic structures, see \cite{IRTU01}. Note, that these manifolds are far from
being spherically generated, for their universal coverings are contractible.
However, some of them satisfy the condition of Theorem~1.

Suppose $M\subseteq \C P^N$ is a symplectic submanifold and let $X$ denote
the blowup of $\C P^N$ along $M$. In \cite{H} the
$\boe$--module structure of $H^*(X)$ is explicitly computed in terms
of the $\boe$--module structure of $H^*(M)$. More precisely, as
$\boe$--modules we have
$$
H^*(X)=H^*(\C P^N)\oplus(H^*(M)\otimes W),
$$
where $W$ is the $\boe$--module $H^*(\C P^{k-2})[2]$ and $2k$ the 
codimension of $M$ in $\C P^N$. For the filtration we therefore get:
$$
H^*(X)_m=
\begin{cases}
H^*(\C P^N)\oplus(H^*(M)_{N-k}\otimes W) & m=N \\
H^*(M)_{m-k}\otimes W & m\neq N 
\end{cases}
$$
It follows from this computation
that $X$ satisfies the Hard Lefschetz Theorem iff $M$ does, \cf
\cite{M84}. As another consequence of this computation we obtain, that the
blow up $X$ satisfies the assumption
of Theorem~\ref{main-theo} iff $M$ does. 

For instance we can take $M$ to be
Thurston's $4$--dimensional nil-manifold embedded in $\C P^5$.
This was the first example of a symplectic manifold which does not satisfy
the Hard Lefschetz Theorem and thus can't be K\"ahler, see \cite{T76}.
The blowup $X$ of $\C P^5$ along $M$ then does not satisfy the Hard 
Lefschetz Theorem for $M$ doesn't.
This was the first example of a closed simply connected
symplectic manifold, which does not satisfy the Hard Lefschetz
Theorem, see \cite{M84}.
However it satisfies the condition of Theorem~\ref{main-theo}, for $M$ does.
So every Hamiltonian fibration with fiber $X$ c--splits. We thus recover
K\c edra's theorem, see \cite{Ka}, without using additional analytic tools.


\begin{thebibliography}{IRTU01}

\bibitem[AB84]{AB84}
M. F. Atiyah and R. Bott,
{\it The moment map and equivariant cohomology},
Topology {\bf 23} (1984) 1--28.

\bibitem[B56]{B56}
A. Blanchard,
{\it Sur les vari\'et\'es analytiques complexes},
Ann. Sci. Ecole Norm. Sup. {\bf 73} (1956), 157--202.

\bibitem[B88]{B88}
J.--L. Brylinski,
{\it A differential complex for {Poisson} manifolds},
J. Diff. Geom. {\bf 28} (1988), 93--114.

\bibitem[H]{H}
S. Haller,
{\it Harmonic cohomology of symplectic manifolds},
to appear in Adv. Math.

\bibitem[IRTU01]{IRTU01}
R. Ib\'a\~nez, Yu. Rudyak, A. Tralle and L. Ugarte,
{\it On symplectically harmonic forms on six-dimensional nil-manifolds},
Comment. Math. Helv. {\bf 76} (2001), 89--109.

\bibitem[Ka]{Ka}
Jaros\l aw K\c edra,
{\it Restrictions on symplectic fibrations},
preprint {\tt math.SG/0203232}.

\bibitem[Kb]{Kb}
Jaros\l aw K\c edra,
{\it Evaluation and universal fibrations},
preprint.

\bibitem[LM02]{LM02}
F. Lalonde and D. McDuff
{\it Symplectic structures on fiber bundles},
to appear in Topology {\bf 42} (2002),  309--347, preprint 
{\tt math.SG/0010275}.

\bibitem[LMP99]{LMP99}
F. Lalonde, D. McDuff and L. Polterovich,
{\it Topological rigidity of Hamiltonian loops and quantum homology},
Invent. Math. {\bf 135} (1999), 369--385.

\bibitem[M95]{M95}
O. Mathieu,
{\it Harmonic cohomology classes of symplectic manifolds},
Comment. Math. Helv. {\bf 70} (1995), 1--9.

\bibitem[M84]{M84}
D. McDuff,
{\it Examples of simply-connected symplectic non-K\"ahlerian manifolds},
J. Diff. Geom. {\bf 20} (1984), 267--277.

\bibitem[M00]{M00}
D. McDuff,
{\it Quantum homology of fibrations over $S^2$},
Inernat. J. Math. {\bf 11} (2000), 665--721.

\bibitem[S01]{S01}
S. Salamon,
{\it Complex structures on nilpotent Lie algebras},
J. Pure Appl. Algebra {\bf 157} (2001), 311--333.

\bibitem[T76]{T76}
W. P. Thurston,
{\it Some simple examples of symplectic manifolds},
Proc. Amer. Math. Soc. {\bf 55} (1976), 467--468.

\end{thebibliography}
\end{document}